\documentclass[12pt]{article}

\makeatletter
\renewcommand{\@oddfoot}{\hfill \thepage}
\makeatother

\usepackage[english]{babel}
\usepackage{mathrsfs}
\usepackage{amssymb}
\usepackage{amsmath,amsthm,amsfonts}
\usepackage{cite}
\usepackage[dvipdf]{graphicx}
\usepackage{fancyhdr}
\usepackage[titletoc]{appendix} 
\usepackage{epstopdf}

\begin{document}

\begin{center}
\Large{\bf Marginals of the planar symmetric Markov random flight on long time intervals
behave like \\ the Goldstein-Kac telegraph process} 
\end{center}

\begin{center}
Alexander D. KOLESNIK\\
Institute of Mathematics and Computer Science\\
Moldova State University\\
Academy Street 5, Kishinev 2028, Moldova\\
\end{center}

\begin{abstract}
The planar symmetric Markov random flight $\bold X(t), \; t>0,$ is represented by the stochastic motion of a particle moving with constant finite speed $c>0$ in the Euclidean plane $\Bbb R^2$ and taking on its initial and each new directions at $\lambda$-Poisson ($\lambda>0$) distributed random time instants by choosing them at random according to the uniform distribution on the unit circumference. We consider the marginals of 
$\bold X(t)$, that is, the projection of this stochastic motion onto the axes. This projection onto the 
$x_1$-axis (respectively, onto the $x_2$-axis) represents a one-dimensional stochastic motion with random velocity $c \cos\alpha$ (respectively, with random velocity $c \sin\alpha$), where $\alpha$ is a random variable distributed uniformly on the interval $[0, 2\pi)$. We prove that the density of the marginals of 
$\bold X(t)$ is asymptotically, as $t\to\infty$, equivalent to the density of the classical one-dimensional Goldstein-Kac telegraph process with parameters ($c, \lambda$). This unexpected and interesting result is confirmed by numerical calculations. 
\end{abstract}

\vskip 0.1cm


\bigskip 


\section{Introduction} 

The theory of Markov random flights in the Euclidean spaces is one of the actively developing areas of stochastic processes in recent decades. Historically, this field of stochastic analysis has started from 
the pioneering works by S. Goldstein \cite{gold} and M. Kac \cite{kac}, where a random motion of a particle moving with finite speed on the real line $\Bbb R^1$ and alternating its direction at Poissonian random time instants, was examined. The main result obtained in these works states that the transition density of the motion is the fundamental solution (the Green's function) to the hyperbolic telegraph equation (also called the damped wave equation), and can be found by solving it with respective initial conditions. Such 
stochastic motion, which is now referred to as the Goldstein-Kac telegraph process and representing the 
one-dimensional Markov random flight, as well as its numerous generalizations, have become the subject of 
a lot of researches, both theoretical and applied. Note that in some works, in order to emphasize that the motion is controlled by a Poisson process, it is also called the Poisson-Kac telegraph process.  

The most important feature of the Goldstein-Kac telegraph process is that it can serve as an appropriate model of the one-dimensional transport that generates a diffusion with finite speed of propagation. For this reason, a generalized telegraph process has become a basis for describing many important physical properties of hyperbolic finite-velocity transport \cite{bras,giona1,giona2,giona3,giona4}. Relativistic aspects of the diffusion generated by finite-velocity random motions on the line, as well as their connections with some other physical processes, were considered  in \cite{cane1,cane2,giona}. An approach to modeling the finite-velocity radioactive transfer based on some properties of such random motions was presented in \cite{eon}.  
A model of cosmic microwave background (CMB) radiation described by a telegraph equation with random initial conditions on the surface of the unit sphere in the three-dimensional space, as well as a number of its astrophysical interpretations, were given in recent works \cite{broad1,broad2}. 
At present, the number of works devoted to telegraph processes, telegraph equations and their applications in various fields of science, technology and engineering is really huge and number more than hundred works. The reader interested in the modern theory of telegraph processes and their applications in financial modeling may address to recent monograph \cite{kolrat} and bibliography therein. 

The natural generalization of the one-dimensional Goldstein-Kac telegraph process is the multidimensional Markov random flight represented by the stochastic motion of a particle that moves at finite speed  
(constant norm of the velocity) in the Euclidean space $\Bbb R^m, \; m\ge 2$, and changes its direction at 
Poissonian random time instants by choosing them at random on the unit $(m-1)$-dimensional sphere 
according to the uniform probability law. In studying such processes, the most desirable goal is undoubtedly to obtain their exact distributions. As it turns out, the difficulty of obtaining the distribution of the 
$m$-dimensional Markov random flight depends fundamentally on the dimension of the space $\Bbb R^m$. Exact distributions of such processes were obtained, by different methods, only in the even-dimensional spaces 
$\Bbb R^2$ (see \cite{kol5,kolors,mas,sta2} or \cite[Section 5.2]{kol1}), 
$\Bbb R^4$ (see \cite{kol4} or \cite[Section 7.2]{kol1}) and 
$\Bbb R^6$ (see \cite{kol2} or \cite[Section 8.2]{kol1}). 
Moreover, the distributions of the Markov random flights in the spaces $\Bbb R^2$ and $\Bbb R^4$ are expressed in terms of elementary (exponential) functions. 
As to the most important case of Markov random flight in the three-dimensional Euclidean space 
$\Bbb R^3$ is concerned, its distribution has not yet been obtained, despite several attempts.  

In the last decade, mostly in physical literature, such finite-velocity stochastic processes have been reincarnated as run-and-tumble motions (for more details on this theory and its applications in physics, see, for instance, \cite{angel,dhar,mall,zhang} and bibliographies therein). A quantum counterpart of 
such finite-velocity random motions was considered in \cite{gzyl}. 

In order to study the behaviour of the multidimensional Markov random flights on a unified basis, a general method of integral transforms of their distributions was developed, which works effectively in the Euclidean space $\Bbb R^m$ of any dimension $m\ge 2$ (see \cite{kol3} or \cite[Chapter 4]{kol1}). One of the most important results obtained by means of this method is an explicit formula for the Laplace-Fourier transform of the transition density of the symmetric Markov random flight in arbitrary dimension $m\ge 2$ 
(see \cite[Formula (4.8)]{kol3} or \cite[Formula (4.6.1)]{kol1}). For the particular low dimensions $m=2$ and $m=3$ this explicit general formula gives already known relations \cite[Formulas (12) and (45), respectively]{mas} and \cite[Formula (5.8)]{sta1}. However, inverting this Laplace-Fourier transform, at least in a direct way, and obtaining an explicit form of the transition density in arbitrary dimension 
$m\ge 3$ is a very difficult and, in fact, impracticable problem. Recently one managed to evaluate just the inverse Laplace transform in this explicit formula and to obtain two representations for the characteristic function of the symmetric Markov random flight in arbitrary dimension $m\ge 3$ in the form of two series 
(with respect to Bessel functions and with respect to the powers of time variable), whose coefficients are given by recurrent relations (see \cite{kol0}). 

Despite the large differences in the behaviour of the Markov random flights in the Euclidean spaces of different dimensions, the natural question arises: is there anything common that is inherent in all such processes in different dimensions and somehow connects their basic characteristics? If so, is it possible to somehow express such characteristics of a Markov random flight in the space $\Bbb R^m$ through known characteristics of the Markov random flights in the spaces of adjacent dimensions $\Bbb R^{m+1}$ and 
$\Bbb R^{m-1}$? This question is of a special relevance since, as noted above, the most important case of  Markov random flight in the three-dimensional space $\Bbb R^3$ has not been studied so far, its distribution is not found, while the distributions of Markov random flights in the neighbouring spaces $\Bbb R^2$ and 
$\Bbb R^4$ are well known. This idea has become the basis of a number of works devoted to attempts to approximate the distribution of a three-dimensional Markov random flight by means of the known distributions of Markov random flights in the spaces of dimensions 2 and 4 (see for instance \cite{paass}, bibliography and figures therein). However, the shape of the three-dimensional density obtained by such approximation looks far from that expected. That is why the question concerning the possibility of such an approximation still remains open.

Nevertheless, it seems possible that there may exist relationships between some characteristics of Markov random flights in different dimensions, at least in some low-dimensional Euclidean spaces. In this article we demonstrate the existence of such a connection between the distribution of the symmetric Markov random flight on the real line $\Bbb R^1$ (the Goldstein-Kac telegraph process $X(t), \; t>0$), and the symmetric Markov random flight $\bold X(t), \; t>0,$ in the Euclidean plane $\Bbb R^2$. The main result of the article states that the density of $X(t)$ and the densities of the marginals of $\bold X(t)$ become asymptotically equivalent, as time $t$ increases, and their difference has the order $O(t^{-1}), \; t\to\infty$, at arbitrary fixed point. 

The article is organized as follows. In Section 2, for the reader's convenience, we recall some preliminary known results concerning the transition density of the one-dimensional Goldstein-Kac telegraph process 
$X(t), \; t>0$, as well as the densities of the planar symmetric Markov random flight $\bold X(t), \; t>0$ and its marginals. In Section 3 we prove the asymptotic equivalence, as time $t\to\infty$, of the transition density of $X(t)$ and the densities of the marginals of $\bold X(t)$ and show that the difference between these densities has the order $O(t^{-1}), \; t\to\infty$, at arbitrary fixed point. Section 4 is devoted to the numerical calculations based on the obtained results. A numerical example with calculated values of these densities and their differences, is given. Several figures that clearly demonstrate the proximity of these densities, as time increases, are also presented. Some final remarks and conclusions are given in Section 5.

\section{Preliminaries} 
\setcounter{equation}{0} 

   The classical Goldstein-Kac telegraph process $X(t), \; t>0,$ is represented by the one-dimensional stochastic motion of a particle that, at the initial time instant $t=0$, starts from the origin 
$0\in\Bbb R^1$ of the real line $\Bbb R^1$ in one of two possible (positive or negative) directions and moves with some finite speed $c$. The particle alternates the directions of its motion at random time instants that form a homogeneous Poissonian flow of rate $\lambda>0$. 

At arbitrary time instant $t>0$, the particle, with probability 1, is located in the closed interval 
$[-ct, ct]$, being, with probability $e^{-\lambda t}$, at the terminal points $\pm ct$ of this interval 
and, with probability $1-e^{-\lambda t}$, at its interior $(-ct, ct)$. 

Let $\text{Pr} \{ X(t)<x \}$ be the distribution of $X(t)$ and let 
$f(x,t) = \frac{\partial}{\partial x} \; \text{Pr} \{ X(t)<x\}$ 
be the density of this distribution. The classical result states that density $f(x,t)$ satisfies the 
second-order hyperbolic partial differential equation 
\begin{equation}\label{tel2}
\frac{\partial^2 f(x,t)}{\partial t^2} + 2\lambda \; \frac{\partial f(x,t)}{\partial t} - 
c^2 \; \frac{\partial^2 f(x,t)}{\partial x^2} = 0
\end{equation}
and the initial conditions 
\begin{equation}\label{incond}
f(x,t)\left|_{t=0} = \delta(x), \right. \qquad \left. \frac{\partial f(x,t)}{\partial t}\right|_{t=0} = 0, 
\end{equation} 
where $\delta(x)$ is Dirac delta-function. Equation (\ref{tel2}) is referred to as the {\it telegraph} 
or {\it damped wave equation}. Initial conditions (\ref{incond}) mean that density $f(x,t)$ is the fundamental solution (Green's function) to the telegraph equation (\ref{tel2}). Its explicit form is given by the formula: 
\begin{equation}\label{dens1}
\aligned
f(x,t) & = \frac{e^{-\lambda t}}{2} \left[ \delta(ct+x) + \delta(ct-x) \right]\\
& + \frac{e^{-\lambda t}}{2c} \left[ \lambda I_0\left(
\frac{\lambda}{c} \sqrt{c^2t^2-x^2}\right) + \frac{\partial}{\partial t} \; I_0\left( \frac{\lambda}{c}
\sqrt{c^2t^2-x^2}\right)  \right] \Theta(ct-\vert x\vert) , \\
& \qquad\qquad\qquad x\in (-\infty, \infty), \qquad t>0,
\endaligned
\end{equation}
where 
$I_0(z) = \sum\limits_{k=0}^{\infty} \frac{1}{(k!)^2} \left( \frac{z}{2} \right)^{2k}$
is the modified Bessel function of zero order, $\delta(x)$ is the Dirac delta-function and $\Theta(x)$ is the Heaviside unit-step function. The first term in (\ref{dens1}) represents the singular component of the density concentrated at the terminal points $\pm ct$, while the second term represents its absolutely continuous part concentrated at the open interval $(-ct, ct)$. 

Taking into account that $I_0'(z)=I_1(z)$, density (\ref{dens1}) can be represented in the following alternative form: 
\begin{equation}\label{dens4}
\aligned f(x,t) & = \frac{e^{-\lambda t}}{2} \left[ \delta(ct+x) + \delta(ct-x) \right]\\
& + \frac{\lambda e^{-\lambda t}}{2c} \left[ I_0\left(
\frac{\lambda}{c} \sqrt{c^2t^2-x^2}\right) + \frac{ct}{\sqrt{c^2t^2-x^2}} \; I_1\left( \frac{\lambda}{c}
\sqrt{c^2t^2-x^2}\right)\right] \Theta(ct-\vert x\vert) , \\
& \qquad\qquad\qquad x\in (-\infty, \infty), \qquad t>0,
\endaligned
\end{equation}
where 
$I_1(z) = \sum\limits_{k=0}^{\infty} \frac{1}{k! \; (k+1)!} \left( \frac{z}{2} \right)^{2k+1}$
is the modified Bessel function of first order. 

The reader interested in the modern state of the telegraph processes theory, its generalizations and applications in financial modeling may address to recent monograph \cite{kolrat}. A number of interesting physical aspects of the generalized telegraph processes and their applications to studying the finite-velocity transport phenomena can be found in \cite{giona1,giona2,giona3,giona4}.

The natural two-dimensional counterpart of the Goldstein-Kac telegraph process is represented by the stochastic motion of a particle that, at the initial time instant $t=0$, starts from the origin 
$\bold 0 = (0, 0)\in\Bbb R^2$ of the the Euclidean plane $\Bbb R^2$ and moves with some finite speed $c$ 
(note that $c$ is treated as the constant norm of the velocity). The initial direction is a random 
two-dimensional unit vector uniformly distributed (the Lebesgue probability measure) on the unit 
circumference 
$S_1 = \left\{ \bold x=(x_1, x_2)\in \Bbb R^2: \; \Vert\bold x\Vert^2 = x_1^2+x_2^2=1 \right\}$ .

The particle changes its direction at random time instants that form a homogeneous Poisson 
flow of rate $\lambda>0$. In each of such Poissonian moments the particle instantly takes on 
a new random direction uniformly distributed on $S_1$ independently of its previous direction. 
Each sample path of this motion represents a planar broken line of total length $ct$ composed of 
the segments of $\lambda$-exponentially distributed random lengths and uniformly oriented 
in $\Bbb R^2$. The trajectories of the process are continuous and differentiable almost everywhere. 

Let $\bold X(t)=(X_1(t), X_2(t))$ be the particle's position at 
arbitrary time instant $t>0$. The stochastic process $\bold X(t)$ is referred to as the 
{\it planar symmetric Markov random flight}. 

At arbitrary time instant $t>0$ the particle, with probability 1, is located in the 
circle $\bold B_{ct} = \left\{ \bold x=(x_1, x_2)\in \Bbb R^2 : \;
\Vert\bold x\Vert^2 = x_1^2+x_2^2\le c^2t^2 \right\} $ of radius $ct$.  

Consider the distribution 
$$\text{Pr} \left\{ \bold X(t)\in d\bold x \right\} = \text{Pr}
\left\{ X_1(t)\in dx_1, X_2(t)\in dx_2 \right\} , \qquad \bold x\in\bold B_{ct}, \quad t>0$$ 
of the stochastic process $\bold X(t)$. The most important result states that the density (in the sense 
of generalized functions) of this distribution has the form   
(see \cite[formula (5.2.5)]{kol1}, \cite[formula (21)]{kolors}, \cite{sta2}): 
\begin{equation}\label{3.1.9}
p(\bold x,t) = \frac{e^{-\lambda t}}{2\pi ct} \; \delta(c^2t^2 - \Vert\bold x\Vert^2) + 
\frac{\lambda}{2\pi c} \; \frac{\exp \left( -\lambda t +
\frac{\lambda}{c} \sqrt{c^2t^2 - \Vert\bold x\Vert^2} \right)}{\sqrt{c^2t^2-\Vert\bold x\Vert^2}} \; 
\Theta(ct-\Vert\bold x\Vert) , 
\end{equation}
$$\bold x = (x_1, x_2)\in\Bbb R^2, \qquad \Vert\bold x\Vert^2 = x_1^2+x_2^2, \qquad t>0,$$
where $\delta(x)$ is the Dirac delta-function and $\Theta(x)$ is the Heaviside unit-step function. 

The first term in (\ref{3.1.9}) represents the singular part of density $p(\bold x,t)$ 
and is related to the case when no Poisson events occur in the time interval $(0,t)$ 
and, therefore, the particle does not change its initial direction (the probability of this event is 
$e^{-\lambda t}$).   The singular component is concentrated on the circumference 
$S_{ct} =\partial\bold B_{ct} = \left\{ \bold x=(x_1, x_2)\in \Bbb R^2: \; 
\Vert\bold x\Vert^2 = x_1^2+x_2^2=c^2t^2 \right\}$ of radius $ct$. 

The second term in (\ref{3.1.9}) represents the absolutely continuous part of density $p(\bold x,t)$ and 
is related to the case when at least one Poisson event occurs in the time interval $(0,t)$ and, therefore, 
the particle at least once changes its direction (the probability of this event is $1-e^{-\lambda t}$). 
In this case, the particle, at time $t$, is located strictly in the interior 
$\text{int} \; \bold B_{ct} = \left\{ \bold x=(x_1,x_2)\in \Bbb R^2: \; 
\Vert\bold x\Vert^2 = x_1^2+x_2^2<c^2t^2 \right\}$  
of the circle $\bold B_{ct}$. 

We are interested in the marginals of the planar Markov random flight $\bold X(t)=(X_1(t),X_2(t)), \; t>0$, that is, the projections of this process onto the coordinate axes. It is clear that the projection 
$X_1(t), \; t>0,$ onto $x_1$-axis represents a one-dimensional stochastic motion with random velocity 
$c \cos\alpha$, changing at $\lambda$-Poisson time instants, where $\alpha$ is a random variable distributed uniformly in the interval $[0, 2\pi)$. Similarly, the projection $X_2(t), \; t>0,$ onto $x_2$-axis represents a one-dimensional stochastic 
motion with random velocity $c \sin\alpha$. 

It was proved in \cite[Theorem 5.2.3]{kol1} that both these marginals have the same distribution whose density is given by the formula: 
\begin{equation}\label{4.2.6}
\aligned 
g(x,t) & = \frac{\partial}{\partial x} \; \text{Pr} \{X_1(t)<x\}  
= \frac{\partial}{\partial x} \; \text{Pr} \{X_2(t)<x\} \\
& = \frac{e^{-\lambda t}}{\pi \sqrt{c^2t^2-x^2}} + \frac{\lambda e^{-\lambda t}}{2c} 
\biggl[ I_0\left( \frac{\lambda}{c}\sqrt{c^2t^2-x^2} \right) + 
\bold L_0\left( \frac{\lambda}{c}\sqrt{c^2t^2-x^2} \right) \biggr] ,
\endaligned 
\end{equation}
$$x\in (-ct, \; ct), \qquad t>0,$$
where $I_0(z)$ and $\bold L_0(z)$ are the modified Bessel and Struve functions of order zero, respectively, given by the series
\begin{equation}\label{4.2.7} 
I_0(z) = \sum_{k=0}^{\infty} \frac{1}{(k!)^2} \left( \frac{z}{2} \right)^{2k}, \qquad \bold L_0(z) = 
\sum_{k=0}^{\infty} \frac{1}{\left[\Gamma\left( k+\frac{3}{2} \right)\right]^2} \left(\frac{z}{2} 
\right)^{2k+1} .
\end{equation}
Note that, in contrast with two-dimensional density (\ref{3.1.9}), its marginals (\ref{4.2.6}) 
are absolutely continuous and do not contain any singular component. The first term in (\ref{4.2.6}) 
is the density of being, at time $t$, on the arc of the circumference $S_{ct}$ located to the left of 
point $x$, while the second term is the density of being in the planar area of the open circle 
$\text{int} \; \bold B_{ct}$ located to the left of this point. 
For more details on the planar symmetric Markov flight see \cite[Chapter 5]{kol1}, 
\cite{kolors,mas,sta2}. 

Comparing density (\ref{dens4}) of the Goldstein-Kac telegraph process and density (\ref{4.2.6}) of the marginals of the two-dimensional symmetric Markov flight  $\bold X(t)$, we see that, although these functions are different, nevertheless, they have somewhat resembled forms. First, we note that both 
these densities tend to zero, as time $t$ increases, at any fixed point $x\in (-ct, ct)$. This follows 
from the fact that the exponential function tends to zero faster than modified Bessel and Struve functions tend to infinity, as $t\to\infty$. This fact will be demonstrated in the next section, where asymptotic formulas for modified Bessel and Struve functions with increasing argument will be given. We will show that both these functions have the same asymptotic rate, as time increases and, therefore, they are asymptotically equivalent at arbitrary point $x\in (-ct, ct)$ for large values of time variable $t$.

\section{Asymptotic equivalence of the densities} 

Consider the densities (\ref{dens4}) and (\ref{4.2.6}). In order to make a comparative asymptotic analysis of these densities, we should exclude all the terms that do not have influence on the asymptotics of their difference. First, we note that, since we intend to study the asymptotics at arbitrary interior point  
$x\in (-ct, ct)$, the singular part of density (\ref{dens4}) of the Goldstein-Kac telegraph process vanishes. Second, since both these densities contain the common term 
$I_0\left( \frac{\lambda}{c}\sqrt{c^2t^2-x^2} \right)$ in square brackets and the same factor 
$\lambda e^{-\lambda t}/(2c)$ in front of the brackets, we may exclude this term from the consideration. 

Hence, we need to study the asymptotic behaviour of the functions 
\begin{equation}\label{AsymptoticsTelegraph}
R(x,t) = \frac{\lambda t}{2} \; \frac{e^{-\lambda t}}{\sqrt{c^2t^2-x^2}} \; I_1\left( \frac{\lambda}{c}
\sqrt{c^2t^2-x^2}\right) ,
\end{equation} 
\begin{equation}\label{AsymptoticsMarginals}
Q(x,t) = \frac{e^{-\lambda t}}{\pi \sqrt{c^2t^2-x^2}} + \frac{\lambda \; e^{-\lambda t}}{2c} 
\; \bold L_0\left( \frac{\lambda}{c}\sqrt{c^2t^2-x^2} \right) ,
\end{equation}
as $t\to\infty$, at arbitrary interior point $x\in (-ct,ct)$. Note that, since each of the functions 
(\ref{AsymptoticsTelegraph}) and (\ref{AsymptoticsMarginals}) contains the factor $e^{-\lambda t}$ which tends to zero, as $t\to\infty$, faster than modified Bessel and Struve functions tend to infinity, 
one can conclude that both these functions $R(x,t)$ and $Q(x,t)$  tend to zero, as $t\to\infty$, for arbitrary fixed $x\in (-ct,ct)$ and therefore, their difference tends to zero too. However, we are interested in the rate of this convergence for large values of time variable $t$, that is, in their asymptotics, as $t\to\infty$. 

Our aim is to prove that both these functions 
(\ref{AsymptoticsTelegraph}) and (\ref{AsymptoticsMarginals}) have a similar asymptotic behaviour (with respect to increasing time variable $t$) at arbitrary fixed point $x\in (-ct,ct)$.  
This result is given by the following theorem.

{\bf Theorem 1.} {\it For arbitrary fixed point $x\in (-ct,ct)$ the following asymptotic relations hold:}
\begin{equation}\label{AsymptoticsTelegraph1}
R(x,t) = \frac{1}{2c} \; \sqrt{\frac{\lambda}{2\pi t}} + O(t^{-3/2}) , \qquad t\to\infty ,
\end{equation} 
\begin{equation}\label{AsymptoticsMarginals1}
Q(x,t) = \frac{1}{2c} \; \sqrt{\frac{\lambda}{2\pi t}} + O(t^{-1}) ,  \qquad t\to\infty .
\end{equation}

\begin{proof}
To prove the statement of the theorem, we need the following asymptotic relation for modified Bessel functions (see, for instance, \cite[Formula 8.451(5)]{gr}): 
\begin{equation}\label{AsymptoticsBessel}
I_{\nu}(z) =  \frac{e^z}{\sqrt{2\pi z}} + O(z^{-3/2}) , \qquad |z|\to\infty, \quad \nu\ge 0.
\end{equation} 
Note that the first term of asymptotic formula (\ref{AsymptoticsBessel}) does not depend on  
index $\nu$, while other terms depending on this index are contained in the second term $O(z^{-3/2})$, 
which tends to zero, as $|z|\to\infty$, faster that the first term does so. 

Applying asymptotic formula (\ref{AsymptoticsBessel}) to modified Bessel function in 
(\ref{AsymptoticsTelegraph}), we have: 
$$\aligned 
R(x,t) & = \frac{\lambda t}{2} \; \frac{e^{-\lambda t}}{\sqrt{c^2t^2-x^2}} \; 
I_1\left( \frac{\lambda}{c} \sqrt{c^2t^2-x^2}\right) \\ 
& \sim \frac{\lambda}{2c \; \sqrt{1-\frac{x^2}{c^2t^2}}} \; \frac{1}{\sqrt{2\pi}} \; 
\left( \lambda t \; \sqrt{1-\frac{x^2}{c^2t^2}} \right)^{-1/2}  
\exp\left( -\lambda t + \lambda t \; \sqrt{1-\frac{x^2}{c^2t^2}} \right)  + O(t^{-3/2}) \\ 
& \sim \frac{\lambda}{2c} \; \frac{1}{\sqrt{2\pi}} \; \frac{1}{\sqrt{\lambda t}} + O(t^{-3/2}) \\ 
& = \frac{1}{2c} \; \sqrt{\frac{\lambda}{2\pi t}} + O(t^{-3/2}) ,  \qquad t\to\infty ,
\endaligned$$ 
where we have used the uniform (with respect to $x\in (-ct, ct)$ and $t>0$) estimate 

\begin{equation}\label{UnEst}
\exp\left( -\lambda t + \lambda t \; \sqrt{1-\frac{x^2}{c^2t^2}} \right) < 1. 
\end{equation}
Thus, formula (\ref{AsymptoticsTelegraph1}) is proved. 

Asymptotic relation (\ref{AsymptoticsMarginals1}) can be obtained by several methods, for example, 
by applying an asymptotic formula for unmodified Struve function in terms of Neumann functions 
(see, for instance, \cite[Formula 8.554]{gr}). This way, however, seems fairly complicated due to the necessity of using the connection between modified and unmodified Struve functions and applying 
an asymptotic formula for Neumann functions of imaginary argument.  

Instead, we prefer another, more simple, way of proving (\ref{AsymptoticsMarginals1}) based on 
the relation 
\begin{equation}\label{L0}
\bold L_0(z) = \frac{4}{\pi} \; \sum_{k=0}^{\infty} \frac{(-1)^k}{2k+1} \; I_{2k+1}(z) ,
\end{equation}
expressing the modified Struve function $\bold L_0(z)$ in terms of modified Bessel functions $I_{\nu}(z)$.  

We note that  
\begin{equation}\label{First Term}
\frac{e^{-\lambda t}}{\pi \sqrt{c^2t^2-x^2}} = O(t^{-1}) , \qquad t\to\infty ,
\end{equation}
that is, the first term of function (\ref{AsymptoticsMarginals}) tends to zero, as $t\to\infty$, 
with the rate $t^{-1}$. 

Let us now consider the second term of function (\ref{AsymptoticsMarginals}). 
Applying (\ref{L0}) and taking into account asymptotic relation (\ref{AsymptoticsBessel}), we obtain: 
$$\aligned 
& \frac{\lambda  e^{-\lambda t}}{2c} \; \bold L_0\left( \frac{\lambda}{c}\sqrt{c^2t^2-x^2} \right) \\ 
& = \frac{\lambda  e^{-\lambda t}}{2c} \; \frac{4}{\pi} \; \sum_{k=0}^{\infty} \frac{(-1)^k}{2k+1} \; 
I_{2k+1}\left( \frac{\lambda}{c} \sqrt{c^2t^2-x^2}\right) \\
& \sim \frac{\lambda e^{-\lambda t}}{2c} \; \frac{4}{\pi} \; \frac{1}{\sqrt{2\pi}} \; 
\left( \frac{\lambda}{c} \sqrt{c^2t^2-x^2} \right)^{-1/2} 
\exp\left( \frac{\lambda}{c} \sqrt{c^2t^2-x^2} \right) \; \sum_{k=0}^{\infty} \frac{(-1)^k}{2k+1} 
+ O(t^{-3/2}) \\ 
& \qquad \left( \text{since} \; \sum_{k=0}^{\infty} \frac{(-1)^k}{2k+1} = \frac{\pi}{4} \right) \\
& = \frac{\lambda}{2c} \; \frac{1}{\sqrt{2\pi}} \; 
\left( \lambda t \sqrt{1-\frac{x^2}{c^2t^2}} \right)^{-1/2} 
\exp\left( -\lambda t + \lambda t \sqrt{1-\frac{x^2}{c^2t^2}} \right) + O(t^{-3/2}) \\ 
& \sim \frac{\lambda}{2c} \; \frac{1}{\sqrt{2\pi}} \; \frac{1}{\sqrt{\lambda t}} + O(t^{-3/2}) \\ 
& = \frac{1}{2c} \; \sqrt{\frac{\lambda}{2\pi t}} + O(t^{-3/2}) ,  \qquad t\to\infty ,
\endaligned$$ 
where we have again used the uniform estimate (\ref{UnEst}) and the asymptotic relation 
$$\left( \lambda t \sqrt{1-\frac{x^2}{c^2t^2}} \right)^{-1/2} \sim (\lambda t)^{-1/2} , 
\qquad t\to\infty,$$ 
valid for arbitrary fixed $x\in (-ct,ct)$.  
Taking into account (\ref{First Term}), we arrive at (\ref{AsymptoticsMarginals1}). 
The theorem is thus completely proved. 
\end{proof}

From Theorem 1 it follows that, at arbitrary fixed point $x\in (-ct, ct)$, the density $f(x,t)$ of 
the Goldstein-Kac telegraph process tends to zero, as $t\to\infty$, faster than the density $g(x,t)$ 
of the marginals of the planar Markov random flight. This leads us to the following main theorem. 

{\bf Theorem 2.}  {\it At arbitrary fixed point $x\in (-ct, ct)$, the following asymptotic formula holds:}
\begin{equation}\label{BasicAsymptotic}
\left| f(x,t) - g(x,t) \right| = O(t^{-1}) , \qquad t\to\infty ,
\end{equation}
{\it where densities $f(x,t)$ and $g(x,t)$ are given by (\ref{dens4}) and (\ref{4.2.6}), respectively.} 

\begin{proof} 
Really, according to Theorem 1, we have: 
$$\aligned 
\left| f(x,t) - g(x,t) \right| & = \left| R(x,t) - Q(x,t) \right| \\
& \sim \left| \frac{1}{2c} \; \sqrt{\frac{\lambda}{2\pi t}} + O(t^{-3/2}) - 
\frac{1}{2c} \; \sqrt{\frac{\lambda}{2\pi t}} - O(t^{-1}) \right| \\ 
& = O(t^{-1}) , \qquad t\to\infty ,  
\endaligned$$
proving (\ref{BasicAsymptotic}). This means that, at arbitrary fixed point $x\in (-ct, ct)$, the estimate 
$$\left| f(x,t) - g(x,t) \right| \le \frac{C}{t} ,$$
holds, where $C$ is some positive constant depending on parameters $c$ and $\lambda$. 
\end{proof}

\section{Numerical examples}

The statement of Theorem 2 can be demonstrated by numerical examples. The results of numerical calculations of the density $f(x,t)$ of the Goldstein-Kac telegraph process and the density $g(x,t)$ of the marginals of the planar Markov random flight, as well as their difference, at fixed point $x=5$ (for the values of parameters $\lambda=1, \; c=2$), are given in Table 1 below.  
\begin{center}
\begin{tabular}{|r||r|r|r|}
\hline
 $t$ & $f(5,t)$ & $g(5,t)$ & $|f(5,t)-g(5,t)|$ \\
\hline\hline
  5 & 0.051002 & 0.050719 & 0.000283 \\
\hline
  10 & 0.046783 & 0.047306 & 0.000523 \\
\hline
  20 & 0.038183 & 0.038526 & 0.000343 \\
\hline
  50 & 0.026466 & 0.026583 & 0.000117 \\
\hline
  100 & 0.019315 & 0.019361 & 0.000046 \\
\hline
  150 & 0.015940 & 0.015965 & 0.000025 \\
\hline
  200 & 0.013878 & 0.013895 & 0.000017 \\
\hline
  300 & 0.011393 & 0.011402 & 0.000009 \\
\hline
  400 & 0.009893 & 0.009899 & 0.000006 \\
\hline
  500 & 0.008863 & 0.008867 & 0.000004 \\ 
\hline
	1000 & 0.006287 & 0.006289 & 0.000002 \\
\hline 
\end{tabular} 
\end{center}
\begin{center}
{\bf Table 1:} {\it Values of densities $f(x,t)$, $g(x,t)$ and their difference at point $x=5$ \\
(for $\lambda=1, \; c=2$)}
\end{center}

We see that, as Theorem 2 states, the difference between these densities (at fixed point $x=5$) 
tends to zero, as time $t$ increases, and the rate of this convergence has the order $t^{-1}$. 

Such behaviour of the difference of these densities are also demonstrated by Figures 1 and 2 below. 

\begin{figure}[htbp]
\begin{minipage}[h]{0.45\linewidth}
\center{\includegraphics[width=1\linewidth]{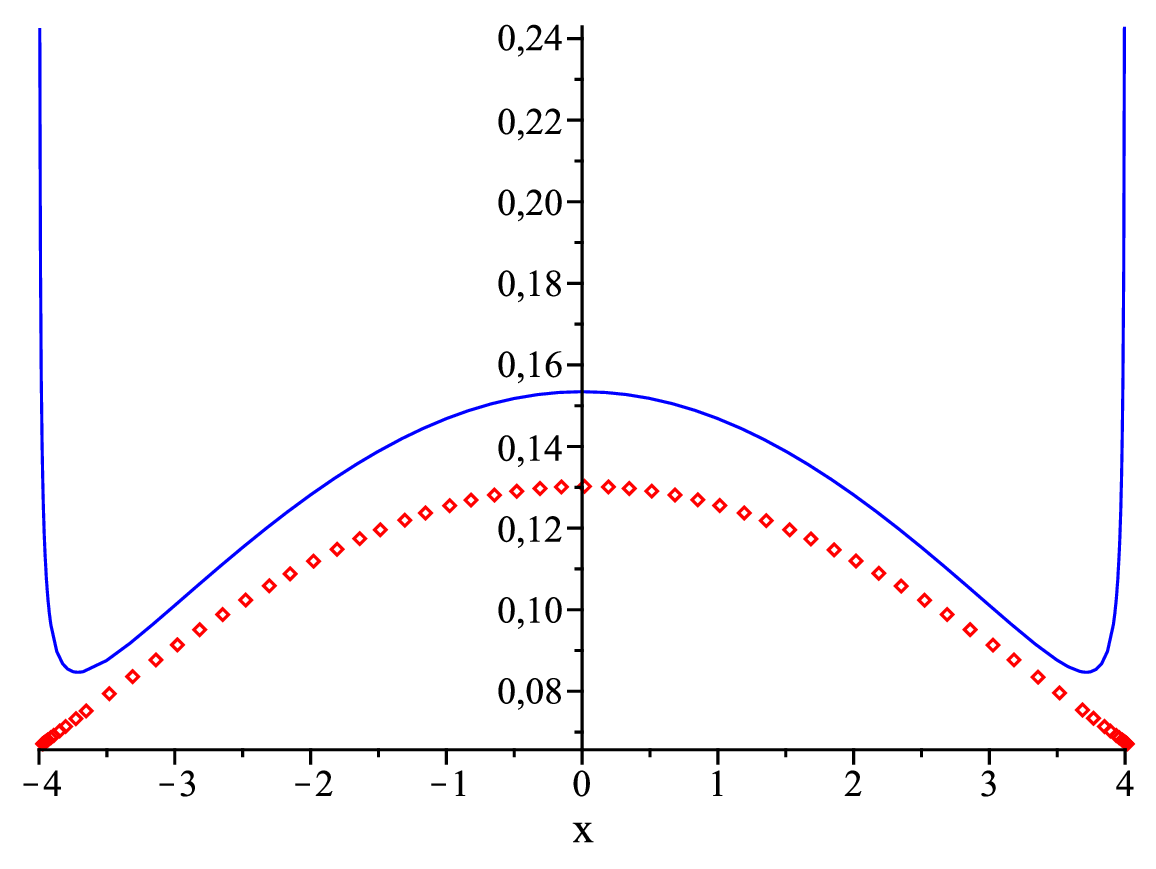}} \\
\end{minipage}
\hfill 
\begin{minipage}[h]{0.45\linewidth}
\center{\includegraphics[width=1\linewidth]{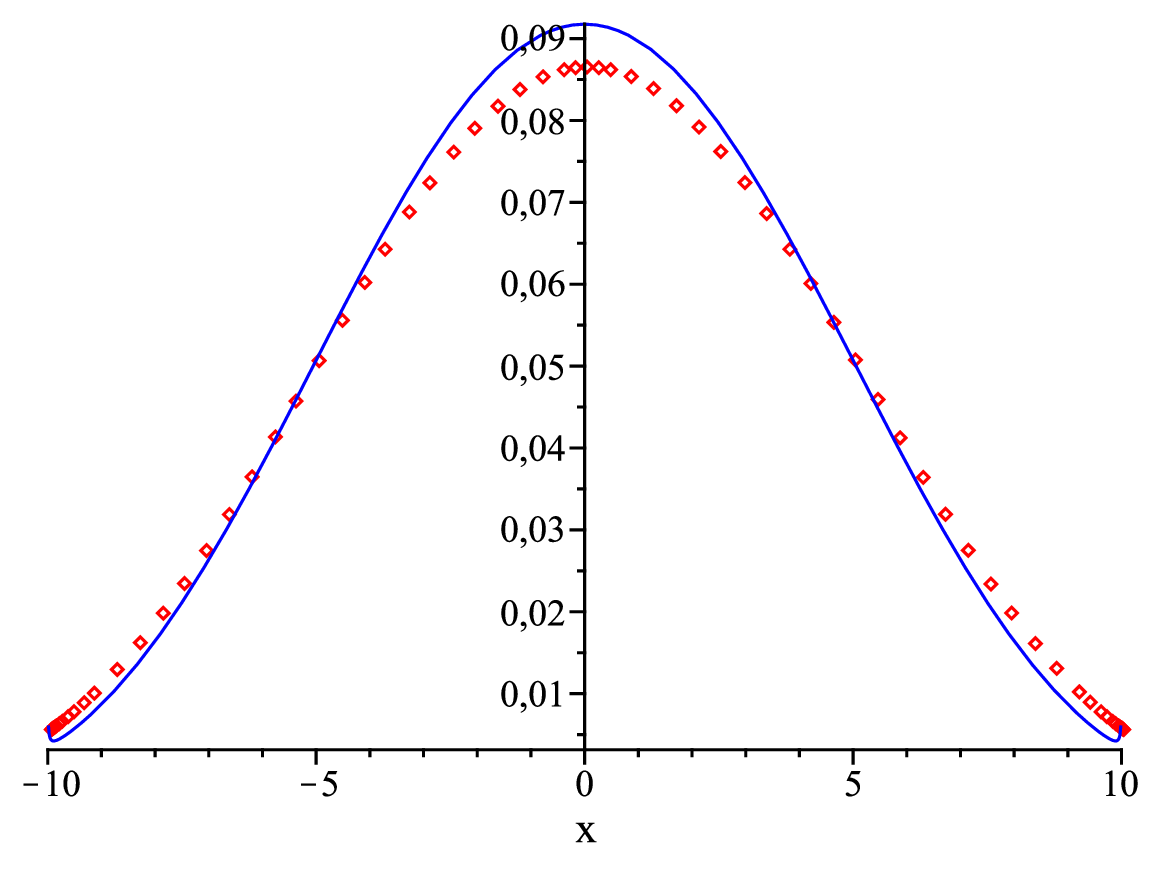}} \\ 
\end{minipage} 
\caption{\it The shapes of densities $f(x,t)$ (dotted line) and $g(x,t)$  
for the time values $t=2$ (left) and $t=5$ (right) (for the values of parameters $\lambda=1, \; c=2$)}   
\label{AsFig1}
\end{figure}
We see that, for the time value $t=2$ (Fig. 1, left), the densities have very different shapes. However, 
for the time value $t=5$ (Fig. 1, right), their shapes become more resembled. 

This resemblance becomes much more obvious for increasing values of time $t$. For the time value 
$t=10$ (Fig. 2, left), the shapes of densities $f(x,t)$ and $g(x,t)$ are very close, while for 
the time value $t=15$ (Fig. 2, right), they are almost identical. It is clear that, for significantly 
larger values of time $t$,  these densities become practically indistinguishable. 

\begin{figure}[htbp]
\begin{minipage}[h]{0.45\linewidth}
\center{\includegraphics[width=1\linewidth]{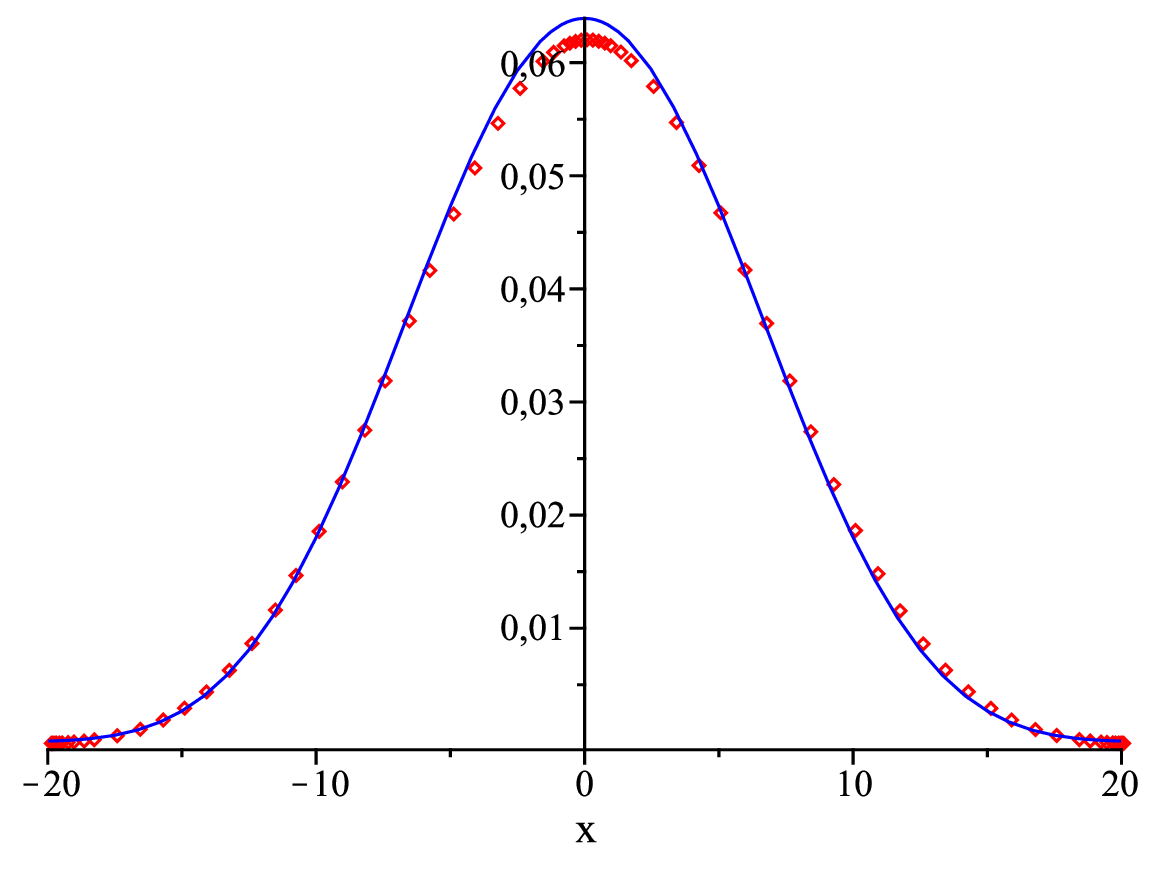}} \\
\end{minipage}
\hfill 
\begin{minipage}[h]{0.45\linewidth}
\center{\includegraphics[width=1\linewidth]{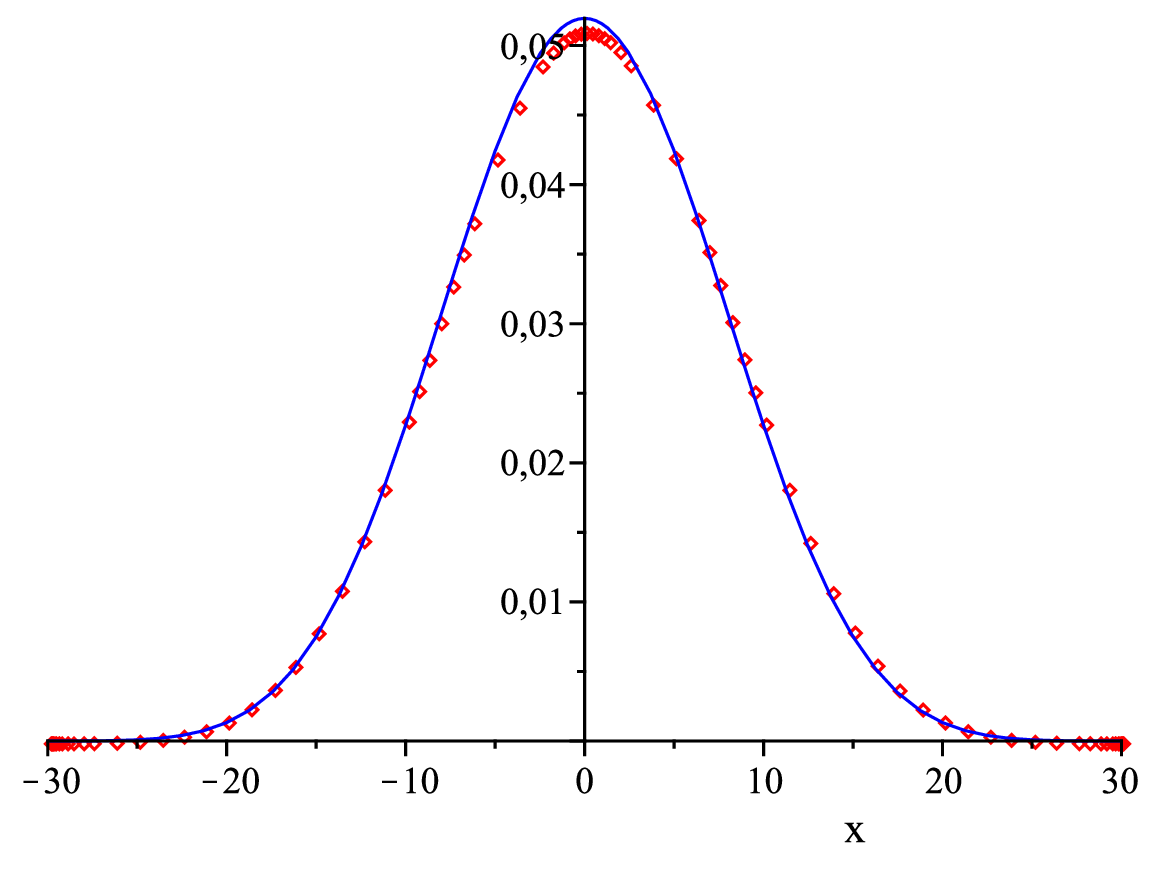}} \\ 
\end{minipage} 
\caption{\it The shapes of densities $f(x,t)$ (dotted line) and $g(x,t)$  
for the time values $t=10$ (left) and $t=15$ (right) (for the values of parameters 
$\lambda=1, \; c=2$)}   
\label{AsFig2}
\end{figure}

The magnitude of their difference, at arbitrary fixed point, is determined by asymptotic formula 
(\ref{BasicAsymptotic}). Note also that these densities have maximal difference at the origin $x=0$.

\section{Conclusions and final remarks}

Analysis of the behaviour of the Markov random flight $\bold X(t), \; t>0,$ in the multidimensional Euclidean spaces $\Bbb R^m, \; m\ge 2$, is of a special importance. It turns out that the complexity 
of such analysis crucially depends on the dimension of the phase space. As it was noted in the 
Introduction, while in the even-dimensional spaces of low dimensions $\Bbb R^2, \; \Bbb R^4$ and 
$\Bbb R^6$ one managed to obtain the distributions of $\bold X(t)$ in explicit forms, the distribution 
of this stochastic process in the three-dimensional Euclidean space $\Bbb R^3$ was not obtained so far, despite many attempts. 
The great importance of studying the three-dimensional Markov random flight on one hand, 
and extreme difficulty in its analyzing, on the other hand, generates the temptation to somehow use 
the known distributions in the neighbouring spaces $\Bbb R^2$ and $\Bbb R^4$ in order to get at least 
some information about the process in the space $\Bbb R^3$. The same concerns all other dimensions. 

The general question may sound as follows: 

\vskip 0.2cm 
{\it Are the distributions of the Markov random flights in neighbouring Euclidean spaces related to each other and, if so, how?} 
\vskip 0.2cm 

In this article we have discovered such connection between the marginals of the planar Markov random 
flight and the Goldstein-Kac telegraph process, namely, that the distributions of these two 
{\bf independent} stochastic processes practically coincide for large values of time variable $t$. 
Moreover, we have also established that the proximity of these distributions has the order $O(t^{-1})$ 
at arbitrary fixed point $x\in (-ct, ct)$. We refer to this property as the {\it asymptotic equivalence} 
of the distributions. 

This amazing and somewhat unexpected result might imply that a similar connection exists in higher dimension, however such a conjecture does not seem to be well founded. The point is that in dimensions 
1 and 2 there are specific properties of some random processes that are no longer observed in higher dimensions. For example, it is known that a Wiener process starting from the origin instantly fills 
both the real line $\Bbb R^1$ and the Euclidean plane $\Bbb R^2$. This means that a Brownian particle, 
at an arbitrarily small moment of time after the start, can, with positive probability, be in  
arbitrarily small interval of the real line $\Bbb R^1$ or in arbitrarily small planar area of the Euclidean plane $\Bbb R^2$, located arbitrarily far from the origin. This phenomenon, however, is not observed in higher dimensions and this fact can be interpreted as the effect of 'density spreading' along additional coordinates.  

Another example, showing a special connection between Markov random flights in the spaces $\Bbb R^1$ (the Goldstein-Kac telegraph process) and $\Bbb R^2$ (the planar Markov random flight), is that the absolutely continuous parts of their densities are the fundamental solutions (the Green's functions) of the identical telegraph equations of respective dimensions (see \cite[Theorems 2.3.1 and 5.4.1]{kol1}). This property is no longer valid for Markov random flights in higher dimensions, whose densities are the fundamental solutions to much more complicated equations than the telegraph one, namely, the hyperparabolic equations whose operators are composed of the multidimensional telegraph and Laplace operators 
(see \cite[Theorem 4.10.1]{kol1}). 

From the results of this article one can conclude that there are some special connections between the symmetric Markov random flights on the real line $\Bbb R^1$ and in the Euclidean plane $\Bbb R^2$. From this fact, however, it does not follow that similar connections exist for Markov random flights in the Euclidean spaces $\Bbb R^m$ of higher dimensions $m\ge3$. 

\bigskip  

{\bf Acknowledgements.} This article was written in the framework of the research project No. 011302.  

\bigskip 

{\bf Declaration.} The author declares no potential conflicts of interest with respect to the research, authorship, and/or publication of this article. The author has no data availability to share.

\end{document}